\documentclass[12pt,letterpaper]{article}
\usepackage{amssymb,amsmath,amsthm}
\usepackage{multirow}
\usepackage{tabularx}
\usepackage[alwaysadjust]{paralist}
\usepackage{filecontents}
\usepackage[utf8x]{inputenc}
\usepackage{subcaption}
\usepackage{graphicx}
\usepackage[numbers,square]{natbib}
\usepackage{tabu}
\usepackage{xcolor}
\usepackage{dcolumn}
\usepackage{sectsty}
\usepackage{booktabs, caption, underscore, url, graphicx}
\sectionfont{\normalsize}
\subsectionfont{\normalsize}
\subsubsectionfont{\normalsize}
\paragraphfont{\normalsize}

\newtheorem*{thm*}{Theorem}

\newcommand{\llll}[1] {\left #1}
\newcommand{\rrrr}[1] {\right #1}

\newcommand{\dddd}[2]{\dfrac{#1}{#2}}

\newcommand{\aaaa}{\alpha}
\newcommand{\tttt}{\tau}

\newcommand{\bbbb}{\beta}
\newcommand{\ssss}{\sigma}

\newcommand{\GGGG}{\Gamma}

\newcommand{\zzzz}{\zeta}
\newcommand{\dddddd}{\delta}

\begin{document}
\nocite{*}
\title{{\bf \normalsize A NEW METHOD FOR NUMERICAL SOLUTION OF THE FRACTIONAL RELAXATION AND SUBDIFFUSION EQUATIONS USING FRACTIONAL TAYLOR POLYNOMIALS}}
\author{Yuri Dimitrov\\
Department of Applied Mathematics and Statistics \\
University of Rousse, Rousse 7017, Bulgaria\\
\texttt{ymdimitrov@uni-ruse.bg}
}
\maketitle
\begin{abstract}
The  accuracy of the numerical solution of a fractional differential equation depends on the  differentiability class of the solution. The derivatives of the solutions of fractional differential equations often have a singularity at the initial point, which may result in a lower accuracy of the numerical solutions. We propose a method for improving the accuracy of the numerical solutions of the fractional relaxation and subdiffusion equations based on the fractional Taylor polynomials
 of the solution at the initial point. 
\end{abstract}
\section{Introduction}

 Differential equations with fractional derivatives are used for modeling non-Markovian processes in various areas of science such as  Physics, Biology and Economics [1-5]. The methods for analytical solution of ordinary and partial fractional differential equations can be applied to only special types of equations and initial conditions. In practical applications the solutions of the fractional differential equations are determined by numerical methods. The development of efficient methods for numerical solution of fractional differential equations is an active research field. 

An important approximation for the Caputo fractional derivative is the $L1$ approximation \eqref{L1}. It has been successfully used for numerical solution of ordinary and partial fractional differential equations. When the function $y(x)$ is a smooth function, the accuracy  of the $L1$ approximation is $O\llll(h^{2-\aaaa}\rrrr)$. In previous work \cite{Dimitrov2015} we determined the second order approximation \eqref{ML1} by modifying the first three coefficients of the $L1$ approximation with the value of the Riemann zeta function at the point $\aaaa-1$, and we computed the numerical solutions of the fractional relaxation and subdiffusion equations with sufficiently differentiable solutions. The solutions of fractional differential equations often have a singularity at the initial point. In this case  the numerical solutions may converge to the exact solution but their  accuracy may be lower than the expected accuracy from the numerical analysis. In this paper we propose a method for improving the accuracy of the numerical solutions of the fractional relaxation equation
\begin{equation}\label{R1}
y^{(\aaaa)}(x)+By(x)=0,\quad y(0)=1,
\end{equation}
where  $0<\aaaa<1$, and the fractional subdiffusion equation
\begin{equation}\label{S1}
	\left\{
	\begin{array}{l l}
\dfrac{\partial^\alpha u(x,t)}{\partial t^\alpha}=\dfrac{\partial^2 u(x,t)}{\partial x^2},\quad x\in [0,\pi],t>0,&  \\
u(x,0)=\sin x,\;u(0,t)=0,\;u(\pi,t)=0.
	\end{array}
		\right . 
	\end{equation}
  When the solutions of the relaxation and subdiffusion equations are smooth functions the accuracy of the numerical solutions using the $L1$ approximation are $O\llll(h^{2-\aaaa}\rrrr)$ and $O\llll(\tttt^{2-\aaaa}+h^2\rrrr)$ and the  numerical solutions using the modified $L1$ approximation have accuracy $O\llll(h^{2}\rrrr)$ and $O\llll(\tttt^{2}+h^2\rrrr)$. The solutions of the fractional relaxation and subdiffusion equations \eqref{R1} and \eqref{S1} have differentiable singularities at the initial point. The accuracy of the numerical solutions \eqref{NS1} and \eqref{NS2} of the fractional relaxation equation is $O\llll(h^{\aaaa}\rrrr)$ and the numerical solutions \eqref{NS3} and \eqref{NS4} of the fractional subdiffusion equation have accuracy $O\llll(\tttt+h^2\rrrr)$. 
	Equations \eqref{R1} and \eqref{S1} are linear fractional differential equations. The form of the equations allows us to compute the Miller-Ross derivatives of the solutions at the initial point. In section 3 and section 4 we use the fractional Taylor polynomials to improve the differentiability properties of the solutions and the accuracy of the numerical methods.
 
The   Caputo  derivative   of order $\aaaa$, where $0<\alpha<1$ is defined as
\begin{equation*}\label{Cder}
D_C^\alpha y(x)=y^{(\alpha)}(x)=\dddd{d^\aaaa}{d x^\aaaa}y(x)=\dfrac{1}{\Gamma (1-\alpha)}\int_0^x \dfrac{y'(\xi)}{(x-\xi)^\alpha}d\xi.
\end{equation*}
When the function $y$ is defined on the interval $(-\infty,x]$, the lower limit of the integral in the definition of Caputo derivative is $-\infty$.  The value of the fractional derivative $y^{(\alpha)}(x)$ depends on the values of the function on the whole interval of definition. One approach for discretizing the Caputo derivative is to divide the interval to subintervals of small length and approximate the derivative of the function $y(x)$ using spline interpolation or a Lagrange polynomial \cite{LiChenYe2011}. Let $h$ be a small positive number. Denote  
$$x_n=n h,\;y_n=y(x_n)=y(n h).$$
 The $L1$ approximation for Caputo derivative is derived \cite{Dimitrov2015} by approximating the derivative of the function $y(x)$ on the interval $\llll[x_{k-1},x_k\rrrr]$ with its value $y^\prime_{k-0.5}$ at the midpoint of the interval.
 \begin{equation} \label{L1}
y^{(\alpha)}_n \approx \dfrac{1}{h^\alpha \GGGG(2-\aaaa)}\sum_{k=0}^{n-1} \ssss_k^{(\alpha)} y_{n-k},
\end{equation}
where $\ssss_0^{(\alpha)}=1,\; \ssss_n^{(\alpha)}=(n-1)^{1-\aaaa}-n^{1-\aaaa}$ and
$$\ssss_k^{(\alpha)}=(k+1)^{1-\alpha}-2k^{1-\alpha}+(k-1)^{1-\alpha},\quad (k=1,\cdots,n-1).$$
When the function $y(x)$ has a continuous second derivative, the accuracy of the $L1$ approximation is $O\llll(h^{2-\aaaa} \rrrr)$ (\cite{LinXu2007}).
The modified $L1$ approximation for the Caputo derivative
\begin{equation}\label{ML1}
y^{(\aaaa)}_n\approx \dddd{1}{\GGGG(2-\aaaa)h^\aaaa}\sum_{k=0}^n \delta_k^{(\aaaa)} y_{n-k},
\end{equation}
has coefficients $\dddddd_k^{(\aaaa)}=\ssss_k^{(\aaaa)}$ for $2\leq k\leq n$,
$$\dddddd_0^{(\aaaa)}=\ssss_0^{(\aaaa)}-\zzzz(\aaaa-1), \dddddd_1^{(\aaaa)}=\ssss_1^{(\aaaa)}+2\zzzz(\aaaa-1),\dddddd_2^{(\aaaa)}=\ssss_2^{(\aaaa)}-\zzzz(\aaaa-1),$$
where $\zzzz(x)$ is the Riemann zeta function. The  modified $L1$ approximation has accuracy  $O\llll(h^{2} \rrrr)$.

By approximating the Caputo derivative at the point $x_n$ with \eqref{L1} and \eqref{ML1} we obtain the numerical solutions \eqref{NS1} and 
\eqref{NS2} of the fractional relaxation equation.
In Table 1 and Table 2 we compute the maximum error and the  order of the numerical solutions \eqref{NS1} and \eqref{NS2} for the following equations 
\begin{equation}\label{R11}
y^{(0.5)}(x)+y(x)=x^2+\dddd{8}{3\sqrt{\pi}}x^{1.5},\quad y(0)=0,
\end{equation}
\begin{equation}\label{R12}
y^{(0.5)}(x)+y(x)=x^{1.25}+\dddd{5\sqrt{2}}{24\pi}\GGGG^2(0.25)x^{0.75},\quad y(0)=0.
\end{equation}
Equations \eqref{R11} and \eqref{R12} have solutions $y(t)=x^2$ and $y(t)= x^{1.25}$. The solution of equation \eqref{R12} has an unbounded second derivative at $t=0$. While the numerical solutions of equation \eqref{R12} converge to the exact solution, their accuracy is smaller than $O\llll(h^{2-\aaaa}\rrrr)$.
			\begin{table}[ht]
    \caption{Maximum error and order of numerical solutions  \eqref{NS1} and \eqref{NS2} for equation \eqref{R11} on the interval $[0,1]$. }
    \begin{subtable}{0.5\linewidth}
      \centering
  \begin{tabular}{l c c }
  \hline \hline
    $h$ & $Error$ & $Order$  \\ 
		\hline \hline
$0.05$         &$0.00281532$     &$1.45619$\\
$0.025$        &$0.00101549$     &$1.47112$\\
$0.0125$       &$0.00036392$     &$1.48049$\\
$0.00625$      &$0.00012986$     &$1.48661$\\
$0.003125$     &$0.00004621$     &$1.49072$\\
\hline
  \end{tabular}
    \end{subtable}
    \begin{subtable}{.5\linewidth}
      \centering
				\quad\quad
  \begin{tabular}{ l  c  c }
    \hline \hline
    $h$ & $Error$ &$Order$  \\ \hline \hline
$0.05$     & $0.000695507$          & $1.90440$ \\
$0.025$    & $0.000182729$          & $1.92836$ \\
$0.0125$   & $0.000047388$          & $1.94711$ \\
$0.00625$  & $0.000012168$         & $1.96139$ \\
$0.003125$ & $3.10\times10^{-6}$    & $1.97206$ \\
\hline
  \end{tabular}
    \end{subtable} 
\end{table}
\vspace{-0.5cm}
\begin{table}[ht]
    \caption{Maximum error and order of numerical solutions  \eqref{NS1} and \eqref{NS2} for equation \eqref{R12} on the interval $[0,1]$. }
    \begin{subtable}{0.5\linewidth}
      \centering
  \begin{tabular}{l c c }
  \hline \hline
    $h$ & $Error$ & $Order$  \\ 
		\hline \hline
$0.05$         &$0.001825790$      &$1.1544$\\
$0.025$        &$0.000806728$      &$1.17836$\\
$0.0125$       &$0.000357787$     &$1.17298$\\
$0.00625$      &$0.000156707$     &$1.19103$\\
$0.003125$     &$0.000067873$     &$1.20717$\\
\hline
  \end{tabular}
    \end{subtable}
    \begin{subtable}{.5\linewidth}
      \centering
				\quad\quad
  \begin{tabular}{ l  c  c }
    \hline \hline
    $h$ & $Error$ &$Order$  \\ \hline \hline
$0.05$     & $0.001825790$           & $1.15440$ \\
$0.025$    & $0.000806728$          & $1.17836$ \\
$0.0125$   & $0.000351853$          & $1.19711$ \\
$0.00625$  & $0.000151948$          & $1.21139$ \\
$0.003125$ & $0.000065135$         & $1.22206$ \\
\hline
  \end{tabular}
    \end{subtable} 
\end{table}

The fractional relaxation equation \eqref{R1}  and  the fractional subdiffusion equation \eqref{S1} have solutions
$$y(x)=E_\aaaa \llll(-Bx^\aaaa\rrrr),\quad u(x,t)=\sin  x E_\aaaa\llll(-t^\aaaa\rrrr).$$
 Both solutions have a differentiable singularity at the initial point of fractional differentiation.
Jin, Lazarov and Zhou \cite{Lazarov20151} discuss the numerical solution of the fractional subdiffusion equation with accuracy in the time direction  $O(\tttt)$ and determine that the accuracy of the numerical solution \eqref{NS1}  of the relaxation equation \eqref{R1} is $O\llll(h^\aaaa\rrrr)$. 
 Murillo and Yuste \cite{MurilloYuste2013} use a finite difference method with adaptive non-uniform timesteps for numerical solution of the fractional subdiffusion and diffusion-wave equations. 
In section 3 and section 4 we present a method for improving the accuracy of the numerical solutions of equations \eqref{R1} and \eqref{S1} based on the fractional Taylor polynomials of the solution.

\section{Preliminaries}
The values of a differentiable function close to the point $x=x_0$ are approximated by its Taylor polynomials of degree $m$
$$y(x)\approx \sum_{n=0}^m \dddd{(x-x_0)^n}{n!}y^{(n)}(x_0).$$
The  definition of fractional Taylor polynomials at the point $x=0$ is based on the composition of Caputo derivatives.
The Miller-Ross sequential fractional derivative for the  Caputo derivative is defined as
$$D^{\aaaa_1}y(x)=D_C^{\aaaa_1}y(x),\quad D^{\aaaa_1+\aaaa_2}y(x)=D_C^{\aaaa_1}D_C^{\aaaa_2}y(x),$$
and
$$D^{\aaaa_1+\aaaa_2+\cdots+\aaaa_n}y(x)=D_C^{\aaaa_1}D_C^{\aaaa_2}\cdots D_C^{\aaaa_n}y(x).$$
When $y(x)$ is a differentiable function
$$D^{0.5+0.5}y(x)=D_C^{0.5}y(x)D_C^{0.5}y(x)=D_C^{1}y(x)=y^\prime(x).$$
The Caputo derivative $D_C^{1}y(x)=y^\prime(x)$ and the Miller-Ross derivative $D^{0.5+0.5}y(x)$ of the function $y(x)=1+x^{0.5}+x$ are different.  
$$D_C^{1}y(x)=y^\prime(x)=1+0.5 x^{-0.5},$$
$$D^{0.5+0.5}y(x)=D_C^{0.5}D_C^{0.5}y(x)=
D_C^{0.5}\llll(\GGGG(1.5)+\dddd{1}{\GGGG(1.5)}x^{0.5}\rrrr)=1.$$
Two important functions in Fractional Calculus are the Mittag-Leffler and the Gamma function
$$\GGGG(x)=\int_0^\infty e^{-\xi}\xi^{x-1}d\xi.$$
The Gamma function has values $\GGGG(0)=\GGGG(1)=1,\GGGG(0.5)=\sqrt \pi$ and satisfies
$$\GGGG(1+x)=x\GGGG(x),\quad \GGGG(x)\GGGG(1-x)=\dddd{\pi}{\sin\pi x}.$$
The one-parameter and two-parameter Mittag-Leffler functions are defined for $\aaaa>0$ as
$$E_\alpha (x)=\sum_{n=0}^\infty \dfrac{x^n}{\Gamma(\alpha n+1)}, \quad
E_{\alpha,\beta} (x)=\sum_{n=0}^\infty \dfrac{x^n}{\Gamma(\alpha n+\beta)}.$$
When $\aaaa=\bbbb=1$ the Mittag-Leffler functions are equal to the exponential function.
The  Mittag-Leffler functions appear in the solutions of fractional differential equations.  The  Laplace transform of the derivatives of the Mittag-Leffler functions satisfy
$$\mathcal{L}\{ x^{\aaaa k+\bbbb-1} E^{(k)}_{\aaaa,\bbbb}(\pm a x^\aaaa) \}(s)=
\dddd{k! s^{\aaaa-\bbbb}}{\llll( s^\aaaa \mp a\rrrr)^{k+1}}.$$
 The next theorem is a generalization of the Mean-Value Theorem for differentiable functions.
\begin{thm*} (Generalized Mean-Value Theorem)
$$y(x)=y(0)+\dddd{x^\aaaa}{\Gamma (\aaaa+1)} y^{(\aaaa)}(\xi_x)\quad (0\leq\xi_x\leq x).$$
\end{thm*}
In this paper we  denote
 $$D^{n\aaaa}y(x)=D^{\aaaa+\aaaa+\cdots+\aaaa}y(x)=D_C^{\aaaa}D_C^{\aaaa}\cdots D_C^{\aaaa}y(x).$$
When the function $y(x)$ has nonzero Miller-Ross derivatives $D^{n\aaaa}y(0)$ the value of $y(h)$ is approximated by the fractional Taylor polynomials \cite{AjouArqubZhourMomani2013}.
$$y(h)\approx \sum_{n=0}^m \dddd{h^{\aaaa n}}{\GGGG(\aaaa n+1)}
D^{n\aaaa}y(0),$$
where $h$ is a small positive number. The
Miller-Ross derivatives of the function $E_\aaaa\llll(x^\aaaa\rrrr)$ satisfy
$$D^{n\aaaa}E_\aaaa\llll(x^\aaaa\rrrr)=E_\aaaa\llll(x^\aaaa\rrrr), \quad \llll. D^{n\aaaa}E_\aaaa\llll(x^\aaaa\rrrr)\rrrr|_{x=0}
=E_\aaaa\llll(0\rrrr)=1.$$
The integer order derivatives of the function $E_\aaaa\llll(x^\aaaa\rrrr)\rightarrow \pm \infty$ as $x\rightarrow 0$.
When $h$ is a small number, the value of $E_\aaaa\llll(h^\aaaa\rrrr)$ is approximated by its fractional Taylor polynomials. 

The fractional relaxation and subdiffusion equations,
\begin{equation}\label{R2}
y^{(0.5)}(x)+y(x)=0,\quad y(0)=1,
\end{equation}
	\begin{equation}\label{S2}
	\left\{
	\begin{array}{l l}
	\dfrac{\partial^{0.5} u(x,t)}{\partial t^{0.5}}=\dfrac{\partial^2 u(x,t)}{\partial x^2},&  \\
	u(x,0)=\sin x,\;u(0,t)=u(\pi,t)=0,&  \\
	\end{array} 
		\right . 
	\end{equation}
	have solutions $y(x)=E_{0.5} \llll(-x^{0.5}\rrrr)$ and $u(x,t)=\sin  x E_{0.5}\llll(-t^{0.5}\rrrr)$. 
	\vspace{-0.5cm}
	\begin{figure}[ht]
\centering
\begin{minipage}{.47\textwidth}
\vspace{0.1cm}
  \centering
  \includegraphics[width=.8\linewidth]{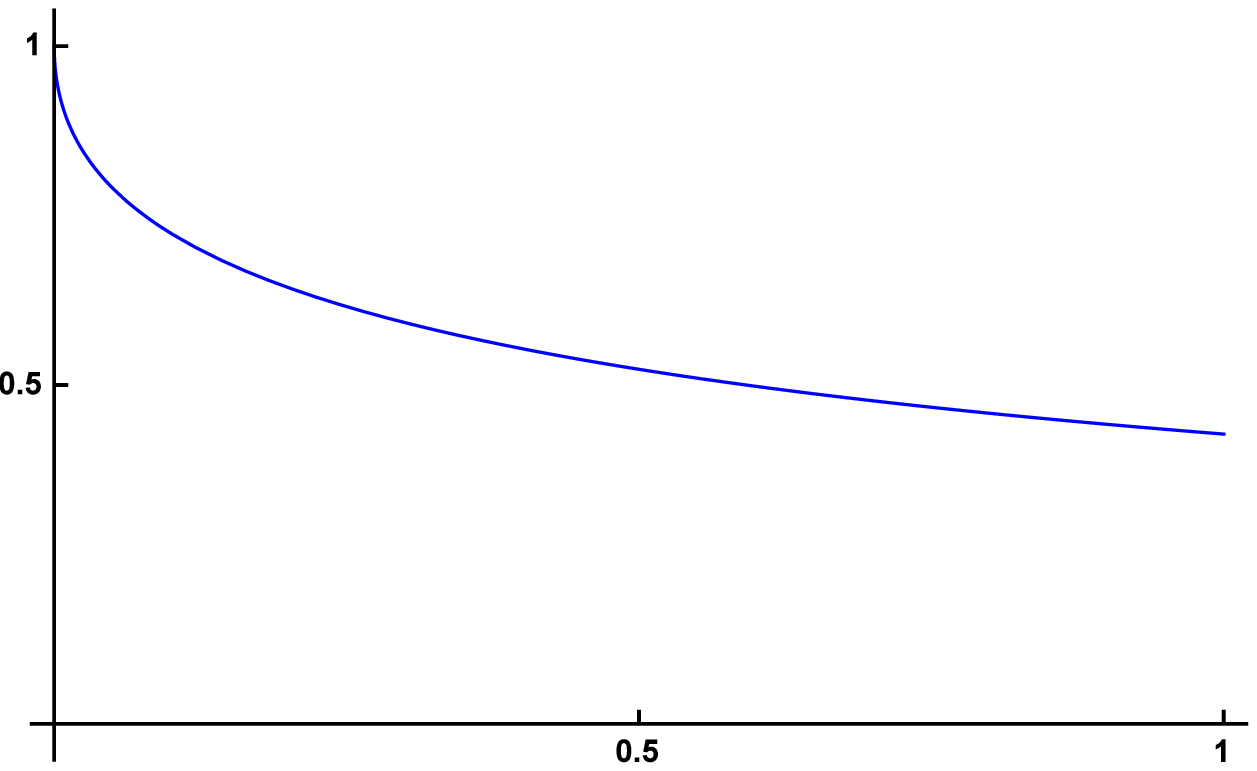}
  \label{fig:test3}
\end{minipage}%
\hspace{0cm}
\begin{minipage}{.47\textwidth}
  \centering
  \includegraphics[width=.8\linewidth]{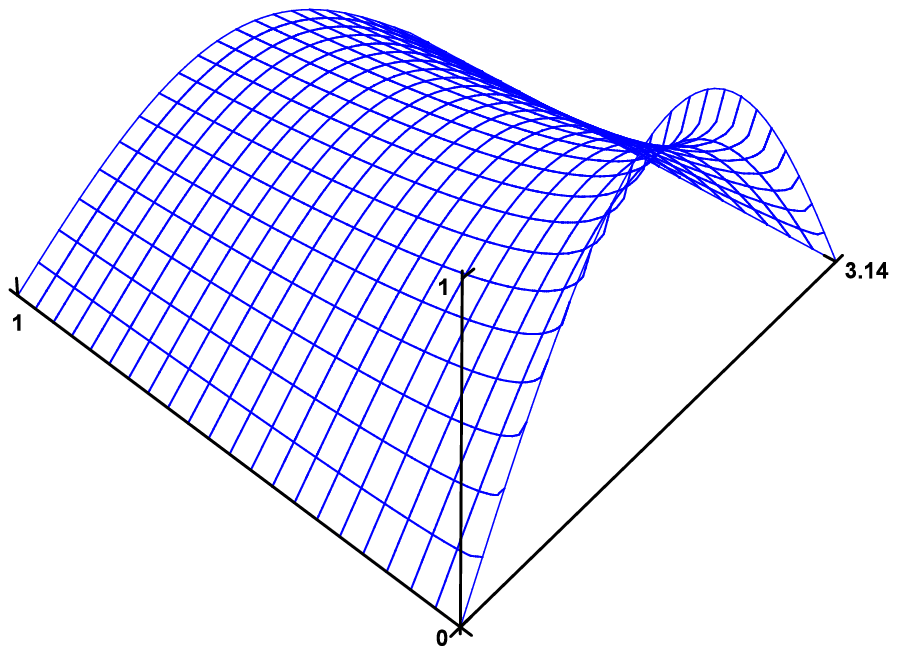}
  \label{fig:test4}
\end{minipage}
 \captionof{figure}{Graphs of the solutions of the fractional relaxation and subdiffusion equations \eqref{R2} and \eqref{S2}.}
\end{figure}
\vspace{-0.6cm}

In \cite{Dimitrov2015} we determined the finite-difference approximations \eqref{NS3} and \eqref{NS4} for the fractional subdiffusion equation. The second derivative in the space direction is approximated by second order central difference. Numerical solution \eqref{NS3} uses the $L1$ approximation for the Caputo derivative in the time direction and \eqref{NS4} uses the modified $L1$ approximation. In Table 3 and Table 4 we compute the maximal error and the order of  numerical solutions \eqref{NS1} and \eqref{NS2} of the fractional relaxation equation \eqref{R2} on the interval $[0,1]$, and the numerical solutions \eqref{NS3} and \eqref{NS4} of the subdiffusion equation \eqref{S2} on the rectangle  $[0,1]\times [0,\pi]$.
	\begin{table}[ht]
    \caption{Maximum error and order of numerical solutions \eqref{NS1} and \eqref{NS3} of the relaxation equation \eqref{R2} on $[0,1]$ (left) and  the subdiffusion equation \eqref{S2} when $\tttt=3h/\pi$, at time $t=1$ (right). }
    \begin{subtable}{0.5\linewidth}
      \centering
  \begin{tabular}{l c c }
  \hline \hline
    $h$ & $Error$ & $Order$  \\ 
		\hline \hline
$0.05$         &$0.0442319$      &$0.378959$\\
$0.025$        &$0.0331978$      &$0.414001$\\
$0.0125$       &$0.024484$     &$0.439247$\\
$0.00625$      &$0.0178342$     &$0.457192$\\
$0.003125$     &$ 0.0128769$     &$0.469859$\\
		\hline
  \end{tabular}
    \end{subtable}
    \begin{subtable}{.5\linewidth}
      \centering
				\quad\quad
  \begin{tabular}{ l  c  c }
    \hline \hline
    $\tttt\;(h=\pi \tttt/3)$ & $Error$ &$Order$  \\ \hline \hline
$0.05\quad$     & $0.00382782$          & $1.08802$ \\
$0.025\quad$    & $0.00184202$          & $1.05524$ \\
$0.0125\quad$   & $0.00089863$          & $1.03549$ \\
$0.00625\quad$  & $0.00044211$          & $1.02333$ \\
$0.003125\quad$ & $0.00021867$          & $1.01563$ \\
\hline
  \end{tabular}
    \end{subtable} 
\end{table}
	\begin{table}[ht]
   \caption{Maximum error and order of numerical solutions \eqref{NS2} and \eqref{NS4} of the relaxation equation \eqref{R2} on $[0,1]$ (left) and  the subdiffusion equation \eqref{S2} when $\tttt=3h/\pi$, at time $t=1$ (right). }
    \begin{subtable}{0.5\linewidth}
      \centering
  \begin{tabular}{l c c }
  \hline \hline
    $h$ & $Error$ & $Order$  \\ 
		\hline \hline
$0.05$         &$0.0442319$      &$0.378959$\\
$0.025$        &$0.0331978$      &$0.414001$\\
$0.0125$       &$0.024484$     &$0.439247$\\
$0.00625$      &$0.0178342$     &$0.457192$\\
$0.003125$     &$ 0.0128769$     &$0.469859$\\
		\hline
  \end{tabular}
    \end{subtable}
    \begin{subtable}{.5\linewidth}
      \centering
				\quad\quad
  \begin{tabular}{ l  c  c }
    \hline \hline
    $\tttt\;(h=\pi \tttt/3)$ & $Error$ &$Order$  \\ \hline \hline
$0.05\quad$     & $0.00222875$          & $1.03303$ \\
$0.025\quad$    & $0.00108596$           & $1.03727$ \\
$0.0125\quad$   & $0.00053127$         & $1.03146$ \\
$0.00625\quad$  & $0.00026116$          & $1.02448$ \\
$0.003125\quad$ & $0.00012893$          & $1.01838$ \\
\hline
  \end{tabular}
    \end{subtable} 
\end{table}

The accuracy of numerical solution \eqref{NS1}  is $O\llll( h^{\aaaa}\rrrr)$ ([16, Example 3.1]).
The maximum error and order of numerical solutions \eqref{NS1} and \eqref{NS2} for the fractional relaxation equation \eqref{R2} are given in Table 3 (left) and Table 4 (left). The two numerical solutions have the same maximal error. The maximum error is attained at the first approximation for $y(h)$. Both numerical solutions use the same approximation $v_1=w_1$ for $y(h)$. Now we show that the error of the approximation $v_1$ for $y(h)$ is $O\llll( h^{0.5}\rrrr)$.
$$y(h)=E_{0.5}\llll(- h^{0.5}\rrrr)=\sum_{n=0}^\infty \dddd{\llll(-h^{0.5}\rrrr)^n}{\GGGG(0.5n+1)},$$
$$y(h)=1-\dddd{2}{\sqrt \pi}h^{0.5}+O\llll(h\rrrr).$$
The approximation $v_1=w_1$ for $y(h)$ is computed with $v_0=1$ and
$$\dddd{v_1-v_0}{\GGGG(1.5)h^{0.5}}+ v_1=0, \quad v_1\llll(1+\GGGG(1.5)h^{0.5}\rrrr)-1=0,$$
$$v_1=\dddd{1}{1+ \GGGG(1.5)h^{0.5}}=1-\dddd{\sqrt \pi}{2}h^{0.5}+O\llll(h\rrrr).$$
Then
$$\llll|y(h)-v_1\rrrr|=\llll(\dddd{\sqrt \pi}{2}-\dddd{2}{\sqrt \pi} \rrrr)h^{0.5}+O\llll(h\rrrr).$$ 
The accuracy of the approximation $v_1=w_1$ for $y(h)$ is $O\llll(h^{0.5} \rrrr)$.
\section{Numerical Solutions of the Fractional Relaxation Equation}
 The fractional relaxation equation  is an important ordinary fractional differential equation. The numerical and analytical solutions of the relaxation equation are discussed in [3,7,9,11-14]. The exact solution of the equation
\begin{equation} \label{RF}
y^{(\aaaa)}(x)+By(x)=F(x),\; y(0)=1,
\end{equation}
 is determined with the Laplace transform method 
$$y(x)= E_{\aaaa}(-B x^\aaaa)+\int_{0}^{x}\xi^{\aaaa-1} E_{\aaaa,\aaaa}\llll(-B \xi^\aaaa\rrrr)F(x-\xi)d\xi.$$
The numerical solutions $\llll\{v_n\rrrr\}$ and $\llll\{w_n\rrrr\}$ of the fractional relaxation equation \eqref{RF} for approximations \eqref{L1} and \eqref{ML1} are computed explicitly with $v_0=1$ (\cite{Dimitrov2015}),
\begin{equation}\label{NS1}
v_n=\dfrac{1}{\ssss_0^{(\alpha)}+B h^\alpha\GGGG(2-\aaaa) } \left( \GGGG(2-\aaaa)h^\alpha F_n-\sum_{k=1}^n \ssss_k^{(\alpha)} v_{n-k} \right),
\end{equation}
and $w_0=1, w_1=v_1$,
\begin{equation}\label{NS2}
w_n=\dfrac{1}{\dddddd_0^{(\alpha)}+B h^\alpha\GGGG(2-\aaaa) } \left( \GGGG(2-\aaaa)h^\alpha F_n-\sum_{k=1}^n \dddddd_k^{(\alpha)} w_{n-k} \right).
\end{equation}
When $F(x)=0$, the fractional relaxation equation 
$$y^{(\aaaa)}(x)+B y(x)=0,\quad y(0)=1,$$
has the solution $y(x)= E_{\aaaa}(-B x^\aaaa)$. 
In section 2 we computed the numerical solutions \eqref{NS1} and \eqref{NS2} of equation \eqref{R1} for $\aaaa=0.5$. Both numerical solutions have maximum error at the first approximation $v_1=w_1$ and accuracy $O\llll( h^{0.5}\rrrr)$. The numerical experiments in this section confirm that the accuracy of $\eqref{NS1}$ and $\eqref{NS2}$  for equation $\eqref{R1}$ is $O\llll( h^{\aaaa}\rrrr)$ when $\aaaa=0.3$ and $\aaaa=0.7$ and all values of $\aaaa$ between $0$ and $1$. In this section we improve the accuracy of $\eqref{NS1}$ and $\eqref{NS2}$ using the fractional Taylor polynomials of the solution at the point $x=0$. The improved numerical solutions have accuracy $O\llll( h^{2-\aaaa}\rrrr)$ and $O\llll( h^{2}\rrrr)$. In Fig. 2 we compare \eqref{NS1} with the improved numerical solution \eqref{INS1} and the exact solution of equation \eqref{R1} for $\aaaa=0.3$.

The fractional relaxation equation is a linear ordinary fractional differential equation. We can determine the Miller-Ross derivatives of the solution at $x=0$ using fractional differentiation. 
The Miller-Ross fractional derivative of order $\aaaa$ is equal to the Caputo derivative. Write the relaxation equation in the form
\begin{equation}\label{EQ}
D^\aaaa y(x)+B y(x)=0,\; y(0)=1.
\end{equation}
When $x=0$ we obtain $D^\aaaa y(0)=-By(0)=-B$. By
applying fractional differentiation of order $\aaaa$
$$D^{2\aaaa} y(x)+BD^{\aaaa}y(x)=0,$$
$$D^{3\aaaa} y(x)+BD^{2\aaaa}y(x)=0.$$
From the above equations we determine the values of $D^{2\aaaa} y(0)$ and $D^{3\aaaa} y(0)$.
$$D^{2\aaaa} y(0)=-BD^{\aaaa}y(0)=B^2,\quad 
D^{3\aaaa} y(0)=-BD^{2\aaaa}y(0)=-B^3.$$
By  differentiating $n-1$ times equation \eqref{EQ} we obtain
$$D^{n\aaaa} y(0)=-BD^{(n-1)\aaaa}y(0).$$
By induction we can show that $D^{n\aaaa} y(0)=(-B)^n$. When $h$ is a small number, the fractional Taylor polynomials approximate the value of $y(h)$
$$y(h)\approx \sum_{n=0}^m \dddd{\llll(-B h^{\aaaa}\rrrr)^n}{\GGGG(\aaaa n+1)}.$$
Now we use the fractional Taylor polynomials for improving the differentiability properties of the solution of the relaxation equation. Let $m$ be a positive integer such that $m\aaaa\geq 2$.  Substitute
$$z(x)=y(x)-\sum_{n=0}^m \dddd{\llll(-B x^{\aaaa}\rrrr)^n}{\GGGG(\aaaa n+1)}.$$
The Caputo derivative of the function $z(x)$ is
$$z^{(\aaaa)}(x)=y^{(\aaaa)}(x)+B\sum_{n=0}^{m-1} \dddd{\llll(-B x^{\aaaa}\rrrr)^n}{\GGGG(\aaaa n+1)}.$$
The function $z(x)$ is a twice continuously differentiable function and satisfies the relaxation equation
\begin{equation}\label{L31}
z^{(\aaaa)}(x)+B z(x)=\dddd{(-B)^{m+1} x^{\aaaa m}}{\GGGG(\aaaa m+1)},\quad z(0)=0.
\end{equation}
The accuracy of numerical solutions \eqref{NS1} and \eqref{NS2} for equation \eqref{L31} is $O\llll( h^{2-\aaaa}\rrrr)$ and $O\llll( h^{2}\rrrr)$. Let $\llll\{\widetilde{z}_k\rrrr\}$ be a numerical solution of \eqref{L31}. The improved numerical solution $\llll\{\widetilde{y}_k\rrrr\}$ of equation \eqref{R1} is determined from $\llll\{\widetilde{z}_k\rrrr\}$ with
\begin{equation}\label{INS1}
\widetilde{y}_k=\widetilde{z}_k+\sum_{n=0}^m (-1)^{n}\dddd{(- B k h)^{0.3 n}}{\GGGG(0.3 n+1)}.
\end{equation}
In Table 5, Table 6 and Table 7 we compute the numerical solutions of the  relaxation equations \eqref{R1} and \eqref{L31} for $\aaaa=0.3$ and $\aaaa=0.7$. The number  $m$ is chosen such that $m\aaaa\geq 2$.

$\bullet$ Let $\aaaa=0.3,B=1$ and $m=7$. The fractional relaxation equation
\begin{equation}\label{L33}
y^{(0.3)}(x)+y(x)=0,\quad y(0)=1
\end{equation}
has analytical solution $E_{0.3}(- x^{0.3})$. Substitute
$$z(x)=y(x)+\sum_{n=0}^7 (-1)^{n+1}\dddd{x^{0.3 n}}{\GGGG(0.3 n+1)}.$$
The function $z(x)\in C^2[0,1]$ and is a solution of the equation
\begin{equation}\label{L34}
z^{(0.3)}(x)+z(x)=\dddd{x^{2.1}}{\GGGG(3.1)},\quad z(0)=0.
\end{equation}
The  accuracy of of numerical solutions \eqref{NS1} and \eqref{NS2} for equation \eqref{L33} is $O(h)$ - Table 5 (left). The maximum error and order of  numerical solutions \eqref{NS1} and \eqref{NS2} for equation \eqref{L34} are given in Table 6. 
In Figure 2 we compare the analytical solution of \eqref{L33} with numerical solution \eqref{NS1} and its improved numerical solution \eqref{INS1}.
\begin{figure}[ht]
  \centering
  \caption{Graph of the analytical solution of equation \eqref{L33}  on the interval $[0,1]$ and  numerical solutions  \eqref{NS1}-red and the improved numerical solution \eqref{INS1}-black for \eqref{NS1}, when $h=0.1$.} 
  \includegraphics[width=0.55\textwidth]{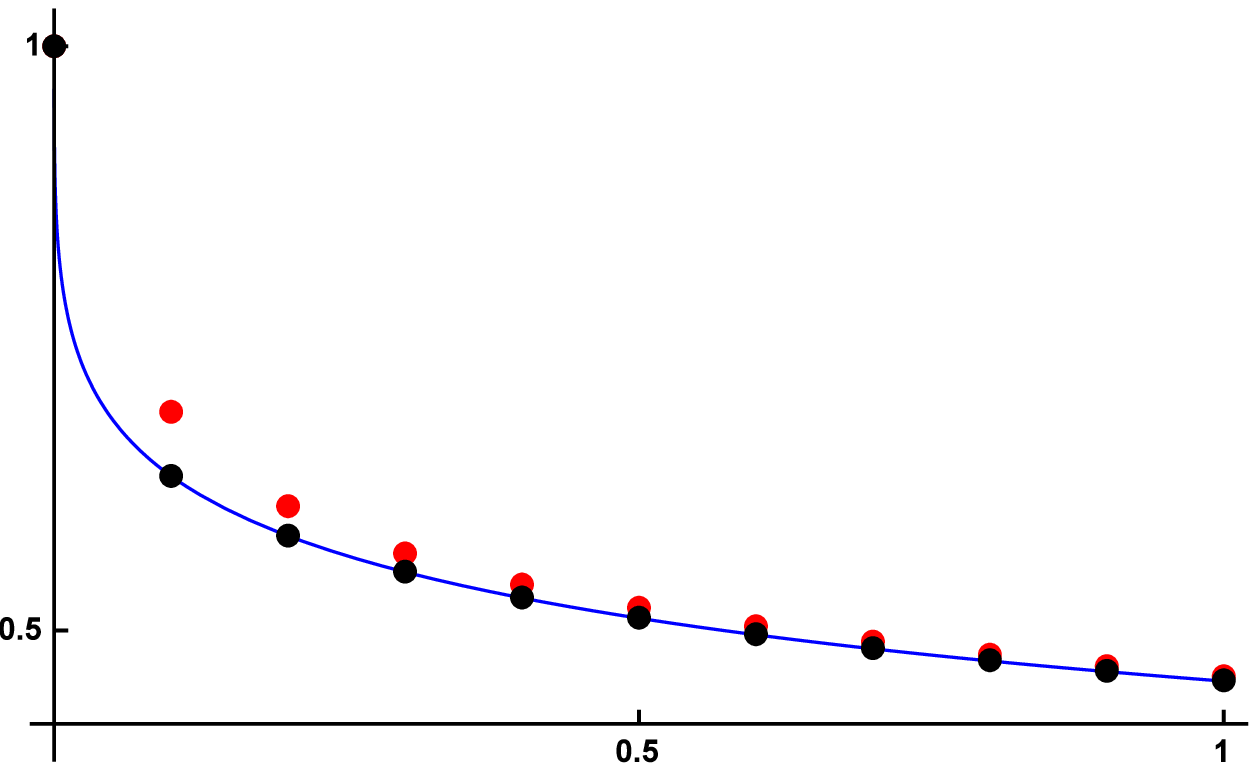}
\end{figure}
	\begin{table}[ht]
    \caption{Maximum error and order of numerical solutions \eqref{NS1} and \eqref{NS2}  for  equation \eqref{L33} (left) and equation \eqref{L35} (right). }
    \begin{subtable}{0.5\linewidth}
      \centering
  \begin{tabular}{l c c }
  \hline \hline
    $h$ & $Error$ & $Order$  \\ 
		\hline \hline
$0.05$         &$0.0496688$     &$0.147528$\\
$0.025$        &$0.0441259$     &$0.170715$\\
$0.0125$       &$0.0386498$     &$0.191164$\\
$0.00625$      &$0.0334386$     &$0.208945$\\
$0.003125$     &$0.0286254$     &$0.224219$\\
		\hline
  \end{tabular}
    \end{subtable}
    \begin{subtable}{.5\linewidth}
      \centering
				\quad\quad
  \begin{tabular}{ l  c  c }
    \hline \hline
    $h$ & $Error$ &$Order$  \\ \hline \hline
$0.05$         &$0.0840166$     &$0.4486824$\\
$0.025$        &$0.0567272$     &$0.566632$\\
$0.0125$       &$0.0365217$     &$0.635288$\\
$0.00625$      &$0.0229408$     &$0.670839$\\
$0.003125$     &$0.0142437$     &$0.687593$\\
\hline
  \end{tabular}
    \end{subtable} 
\end{table}
	
$\bullet$  Let $\aaaa=0.7,m=3$ and $B=4.$
\begin{equation}\label{L35}
y^{(0.7)}(x)+4y(x)=0,\quad y(0)=1.
\end{equation}
The maximum error and order of  numerical solutions \eqref{NS1} and \eqref{NS2} for equation \eqref{L35} are given in Table 5 (right).
Substitute
$$z(x)=y(x)-1+\dddd{4x^{0.7}}{\GGGG(1.7)}-\dddd{16x^{1.4}}{\GGGG(2.4)}+\dddd{64 x^{2.1}}{\GGGG(3.1)}.$$
The function $z(x)$ is a solution of the equation
\begin{equation}\label{L36}
z^{(0.7)}(x)+z(x)=\dddd{256x^{2.1}}{\GGGG(3.1)},\quad z(0)=0.
\end{equation}
The maximum error and order of  numerical solutions \eqref{NS1} and \eqref{NS2} for equation \eqref{L36} are given in Table 7.
\begin{table}[ht]
    \caption{Maximum error and order of numerical solutions \eqref{NS1} and \eqref{NS2}  for  equation \eqref{L33}. }
    \begin{subtable}{0.5\linewidth}
      \centering
  \begin{tabular}{l c c }
  \hline \hline
    $h$ & $Error$ & $Order$  \\ 
		\hline \hline
$0.05$     & $0.00022816$          & $1.61147$ \\
$0.025$    & $0.00007356$          & $1.63301$ \\
$0.0125$   & $0.00002347$          & $1.64829$ \\
$0.00625$  & $7.43\times 10^{-6}$  & $1.65956$ \\
$0.003125$ & $2.34\times 10^{-6}$  & $1.66807$ \\
		\hline
  \end{tabular}
    \end{subtable}
    \begin{subtable}{.5\linewidth}
      \centering
				\quad\quad
  \begin{tabular}{ l  c  c }
    \hline \hline
    $h$    & $Error$               &$Order$  \\ \hline \hline
$0.05$     & $0.0000506703$        & $1.88353$ \\
$0.025$    & $0.0000132354$        & $1.93674$ \\
$0.0125$   & $3.39\times 10^{-6}$  & $1.96378$ \\
$0.00625$  & $8.61\times 10^{-7}$  & $1.97867$ \\
$0.003125$ & $2.17\times 10^{-7}$  & $1.98724$ \\
\hline
  \end{tabular}
    \end{subtable} 
\end{table}

	\begin{table}[ht]
    \caption{Maximum error and order of numerical solutions \eqref{NS1} and \eqref{NS2}  for  equation \eqref{L35}. }
    \begin{subtable}{0.5\linewidth}
      \centering
  \begin{tabular}{l c c }
  \hline \hline
    $h$ & $Error$ & $Order$  \\ 
		\hline \hline
$0.05$     & $0.0825401$          & $1.27437$ \\
$0.025$    & $0.0338612$          & $1.28546$ \\
$0.0125$   & $0.0138328$          & $1.29154$ \\
$0.00625$  & $0.0056374$          & $1.29499$ \\
$0.003125$ & $0.0022943$        & $1.29701$ \\
		\hline
  \end{tabular}
    \end{subtable}
    \begin{subtable}{.5\linewidth}
      \centering
				\quad\quad
  \begin{tabular}{ l  c  c }
    \hline \hline
    $h$ & $Error$ &$Order$  \\ \hline \hline
$0.05$     & $0.00682601$          & $2.21359$ \\
$0.025$    & $0.00136314$          & $2.32411$ \\
$0.0125$   & $0.00024704$          & $2.46412$ \\
$0.00625$  & $0.00004329$          & $2.51272$ \\
$0.003125$ & $7.24\times 10^{-6}$  & $2.57940$ \\
\hline
  \end{tabular}
    \end{subtable} 
\end{table}
\section{Numerical Solutions of the Fractional Subdiffusion Equation}

The time-fractional  subdiffusion equation is a parabolic partial fractional differential equation obtained from the heat equation by replacing the time derivative with a  fractional derivative of order $\aaaa$, where $0<\aaaa<1$,
\begin{equation}\label{L41}
\dfrac{\partial^\alpha u(x,t)}{\partial t^\alpha}=\dfrac{\partial^2 u(x,t)}{\partial x^2}.
\end{equation}
The analytical solution of the fractional subdiffusion equation with initial and boundary conditions
$$u(x,0)=u_0(x),\;u(0,t)=u(\pi,t)=0,$$
is determined using the separation of variables method 
$$u(x,t)=\sum_{n=1}^\infty c_n E_\aaaa (-n^2  t^\aaaa) \sin (n x),$$
where
$$c_n=\dddd{2}{\pi}\int_0^1 u_0(\xi) \sin (n \xi)d\xi.$$
Each term $E_\aaaa (-n^2  t^\aaaa) \sin n x$
 is a solution of \eqref{L41} and the coefficient $c_n$ is the coefficient of the Fourier sine series of the function $u_0(x)$. When $u_0(x)=\sin x$, the subdiffusion equation \eqref{S1} has analytical solution $u(x,t)=\sin  x E_\aaaa (-  t^\aaaa) $. 
In section 3 we used  the fractional Taylor polynomials of the solution of the relaxation equation to improve the accuracy of numerical solutions \eqref{NS1} and \eqref{NS2}. In this section we use the fractional Taylor polynomials of the solution of the subdiffusion equation in the time direction for improving the accuracy of difference approximations \eqref{NS3} and \eqref{NS4}.

Let $h=\pi/N,\tttt=1/M$, where $M$ and $N$ are positive integers, and $\mathcal{G}$ be a grid on the rectangle $[0,\pi]\times[0,1]$
$$\mathcal{G}=\llll\{(n h,m \tttt) \| 1\leq n\leq N, 1\leq m\leq M\rrrr\}.$$
	In \cite{Dimitrov2015} we determined the finite difference approximations \eqref{NS3} and \eqref{NS4} for the  the fractional subdiffusion equation
	\begin{equation*}
	\left\{
	\begin{array}{l l}
	\dfrac{\partial^\alpha u(x,t)}{\partial t^\alpha}=\dfrac{\partial^2 u(x,t)}{\partial x^2}+F(x,t),&  \\
	u(x,0)=u_0(x),\;u(0,t)=u(\pi,t)=0.&  \\
	\end{array} 
		\right . 
	\end{equation*}
Denote by $u_n^m$ and $F_n^m$ the values of the functions $u(x,t)$ and $F(x,t)$ on $\mathcal{G}$.
$$u_n^m=u(n h,m \tttt),\quad F_n^m=F(n h,m\tttt).$$
Let
	$$\eta=\GGGG(2-\aaaa)\dddd{\tttt^\aaaa}{h^2},$$
$K$ be the tridiagonal matrix of dimension $N-1$ with values $ 1+2\eta$ on the main diagonal, and $-\eta$ on the diagonals above and below the main diagonal
		$$K_{5} =
 \begin{pmatrix}
  1+2\eta &  -\eta    & 0        & 0          & 0       \\
 -\eta    &1+2\eta    & -\eta    & 0          & 0       \\
  0       & -\eta     & 1+2\eta  & -\eta      & 0       \\
  0       & 0         & -\eta    &  1+2\eta   & -\eta  \\
	0       & 0         & 0        & -\eta      &  1+2\eta   
 \end{pmatrix}
$$
and $V^0$ be the vector of dimension $N-1$ determined from the initial condition
	$$V^0=\llll[u_0( n h)\rrrr]_{n=1}^{N-1}.$$
	The difference approximation for the subdiffusion equation $\llll\{ V^m\rrrr\}$ using the $L1$ approximation \eqref{L1} is computed with the linear system
	\begin{equation}\label{NS3}
	K V^m=R_1,
	\end{equation}
 where $R_1$ is the vector
	$$R_1=\llll[-\sum_{k=1}^m \ssss_k^{(\aaaa)} V_n^{m-k}+\tttt^\aaaa \GGGG(2-\aaaa)F_n^m \rrrr]_{n=1}^{N-1}.$$
		The numerical solution  $\llll\{W^m\rrrr\}$ of the subdiffusion equation using the modified $L1$ approximation \eqref{ML1} is computed with the linear system
	\begin{equation}\label{NS4}
	(K-\zzzz(\aaaa-1)I) W^m=R_2,
	\end{equation}
  for $2\leq m\leq N-1$, where $R_2$ is the vector
	$$R_2=\llll[-\sum_{k=1}^m \dddddd_k^{(\aaaa)} W_n^{m-k}+\tttt^\aaaa \GGGG(2-\aaaa)F_n^m \rrrr]_{n=1}^{N-1},$$
	and $W^0=V^0,W^1=V^1$.
	
	 When the solution of the subdiffusion equation is a sufficiently differentiable function, the difference approximations \eqref{NS3} and \eqref{NS4} have accuracy $O\llll(\tttt^{2-\aaaa}+h^2 \rrrr)$ and $O\llll(\tttt^{2}+h^2 \rrrr)$. The solution of subdiffusion equation $\eqref{S1}$ has an unbounded first derivative at $t=0$. The numerical experiments in section 2 with $\aaaa=0.5$ and other values of $\aaaa \in (0,1)$ indicate that the accuracy of numerical solutions \eqref{NS3} and \eqref{NS4} is $O\llll(\tttt+h^2 \rrrr)$.
	
Let $D_t^\bbbb u(x,t)$ be the Miller-Ross derivative of order $\bbbb$ of the function $u(x,t)$ in the time direction.
	$$D_t^\aaaa u(x,t)=\dfrac{\partial^2 u(x,t)}{\partial x^2}.$$
	Then
	$$D_t^{2\aaaa}u(x,t)=D_t^{\aaaa}D_t^{\aaaa}u(x,t)=D_t^\aaaa\dfrac{\partial^2 u(x,t)}{\partial x^2}=
	\dfrac{\partial^2 }{\partial x^2}D_t^\aaaa u(x,t)=\dfrac{\partial^4 u(x,t)}{\partial x^4},$$
	$$D_t^{3\aaaa}u(x,t)=D_t^{\aaaa}D_t^{2\aaaa}u(x,t)=D_t^\aaaa\dfrac{\partial^4 u(x,t)}{\partial x^4}=
	\dfrac{\partial^4 }{\partial x^4}D_t^\aaaa u(x,t)=\dfrac{\partial^6 u(x,t)}{\partial x^6}.$$
	In this way we obtain
	$$D_t^{n\aaaa}u(x,t)=\dfrac{\partial^{2n} u(x,t)}{\partial x^{2n}}.$$
Set $t=0$
		$$D_t^{n\aaaa}u(x,0)=\dfrac{\partial^{2n} u(x,0)}{\partial x^{2n}}=
		\dfrac{\partial^{2n} \sin x}{\partial x^{2n}}=(-1)^{n}\sin x.$$
	From the fractional Taylor polynomials approximation
		$$u(x,h)\approx \sin x \sum_{n=0}^m (-1)^n \dddd{h^{n\aaaa}}{\GGGG(n\aaaa+1)},$$
	when $h>0$ is a sufficiently small number. Substitute
		$$v(x,t)=u(x,t)-\sin x \sum_{n=0}^m (-1)^n \dddd{t^{n\aaaa}}{\GGGG(n\aaaa+1)}.$$
		We have that
		$$D_t^\aaaa v(x,t)=D_t^\aaaa u(x,t)+\sin x \sum_{n=0}^{m-1} (-1)^{n+1} \dddd{t^{n\aaaa}}{\GGGG(n\aaaa+1)},$$
		$$\dfrac{\partial^{2}}{\partial x^{2}}v(x,t)=\dfrac{\partial^{2}}{\partial x^{2}}u(x,t)+\sin x \sum_{n=0}^{m} (-1)^{n+1} \dddd{t^{n\aaaa}}{\GGGG(n\aaaa+1)}.$$
	The function $v(x,t)$ is a solution of the fractional subdiffusion equation
		\begin{equation*}
	\left\{
	\begin{array}{l l}
	\dfrac{\partial^{\aaaa} v(x,t)}{\partial t^{\aaaa}}=\dfrac{\partial^2 v(x,t)}{\partial x^2}+(-1)^{m+1}\sin x\dddd{t^{m\aaaa}}{\GGGG(m\aaaa+1)},&  \\
	v(x,0)=v(0,t)=v(\pi,t)=0.&  \\
	\end{array}
		\right . 
	\end{equation*}
When $\aaaa=0.3$ and $m=7$, the fractional subdiffusion equation
	\begin{equation}\label{S03}
	\left\{
	\begin{array}{l l}
	\dfrac{\partial^{0.3} u(x,t)}{\partial t^{0.3}}=\dfrac{\partial^2 u(x,t)}{\partial x^2},&  \\
	u(x,0)=\sin x,\;u(0,t)=u(\pi,t)=0,&  \\
	\end{array} 
		\right . 
	\end{equation}
has analytical solution $u(x,t)=\sin x E_{0.3}\llll(-t^{0.3}\rrrr)$.
The difference approximations \eqref{NS3} and \eqref{NS4} for equation \eqref{S03} have accuracy $O(\tttt)$ in the time direction (Table 8 (left) and Table 9 (left)).
Let
		$$v(x,t)=u(x,t)-\sin x \sum_{n=0}^7 (-1)^n \dddd{h^{0.3n}}{\GGGG(0.3 n+1)}.$$
		The function $v(x,t)$ is a twice continuously differentiable function and is a solution of the homogeneous subdiffusion equation
		\begin{equation}\label{MS03}
	\left\{
	\begin{array}{l l}
	\dfrac{\partial^{0.3} v(x,t)}{\partial t^{0.3}}=\dfrac{\partial^2 v(x,t)}{\partial x^2}+\sin x\dddd{t^{2.1}}{\GGGG(3.1)},&  \\
	v(x,0)=v(0,t)=v(\pi,t)=0.&  \\
	\end{array}
		\right . 
	\end{equation}
In Table 8 (right) and Table 9 (right) we compute the maximum error and numerical order of difference approximations \eqref{NS3} and \eqref{NS4} for the fractional subdiffusion equation \eqref{MS03} on the rectangle 
	$$(x,t)\in [0,\pi]\times[0,1],$$
	with step size in the space direction  $h=\pi \tttt/3$. 
Numerical experiments with other values of the step sizes $h$ and $\tttt$ confirm that the numerical solutions \eqref{NS3} and \eqref{NS4} converge to the analytical solution of the fractional subdiffusion equation \eqref{MS03} with  the expected accuracy $O\llll(\tttt^{2-\aaaa}+h^2 \rrrr)$ and $O\llll(\tttt^{2}+h^2 \rrrr)$. A question for future work is to extend the method presented in this paper to other ordinary and partial fractional differential equations.
		\begin{table}[ht]
    \caption{Maximum error and order of numerical solution \eqref{NS3} of equations \eqref{S03} (left) and  \eqref{MS03} (right), for $h=\pi \tttt/3$, at time $t=1$. }
    \begin{subtable}{0.5\linewidth}
      \centering
  \begin{tabular}{l c c }
  \hline \hline
    $\tttt$ & $Error$ & $Order$  \\ 
		\hline \hline
$0.05$         &$0.00208324$     &$1.08653$\\
$0.025$        &$0.00100782$     &$1.04759$\\
$0.0125$       &$0.00049480$     &$1.02633$\\
$0.00625$      &$0.00024489$     &$1.01468$\\
$0.003125$     &$0.00012175$     &$1.00826$\\
		\hline
  \end{tabular}
    \end{subtable}
    \begin{subtable}{.5\linewidth}
      \centering
				\quad\quad
  \begin{tabular}{ l  c  c }
    \hline \hline
    $\tttt$ & $Error$ &$Order$  \\ \hline \hline
$0.05$     & $0.000247016$          & $1.64590$ \\
$0.025$    & $0.000078268$          & $1.65812$ \\
$0.0125$   & $0.000024643$          & $1.66722$ \\
$0.00625$  & $7.72\times 10^{-6}$   & $1.67411$ \\
$0.003125$ & $2.41\times 10^{-6}$   & $1.67940$ \\
\hline
  \end{tabular}
    \end{subtable} 
\end{table}
		\begin{table}[ht]
    \caption{Maximum error and order of numerical solution \eqref{NS4} of equations \eqref{S03} (left) and  \eqref{MS03} (right), for $h=\pi \tttt/3$, at time $t=1$. }
    \begin{subtable}{0.5\linewidth}
      \centering
  \begin{tabular}{l c c }
  \hline \hline
    $\tttt$ & $Error$ & $Order$  \\ 
		\hline \hline
$0.05$         &$0.001565250$     &$1.05230$\\
$0.025$        &$0.000760933$     &$1.04055$\\
$0.0125$       &$0.000372634$     &$1.03001$\\
$0.00625$      &$0.000183395$     &$1.02281$\\
$0.003125$     &$0.000090563$     &$1.01796$\\
		\hline
  \end{tabular}
    \end{subtable}
    \begin{subtable}{.5\linewidth}
      \centering
				\quad\quad
  \begin{tabular}{ l  c  c }
    \hline \hline
    $\tttt$ & $Error$ &$Order$  \\ \hline \hline
$0.05$     & $0.000031877$        & $1.81175$ \\
$0.025$    & $8.54\times10^{-6}$   & $1.90107$ \\
$0.0125$   & $2.22\times10^{-6}$   & $1.94430$ \\
$0.00625$  & $5.67\times10^{-7}$   & $1.96751$ \\
$0.003125$ & $1.44\times10^{-7}$   & $1.98068$ \\
\hline
  \end{tabular}
    \end{subtable} 
\end{table}

\end{document}